# Karl Weierstrass' Bicentenary

## G.I. Sinkevich

Saint Petersburg State University of Architecture and Civil Engineering Vtoraja Krasnoarmejskaja ul. 4, St. Petersburg, 190005, Russia

galina.sinkevich@gmail.com

*Abstract*. Academic biography of Karl Weierstrass, his basic works, influence of his doctrine on the development of mathematics.

In 2015, the mathematical world celebrates the bicentenary of the great German mathematician Karl Weierstrass (1815-1897) who created the modern mathematical analysis.

***Childhood and adolescence***. Karl Theodor Wilhelm Weierstrass was born on 31 October 1815 in Ostenfelde (Westphalia) to a catholic family of burgomaster's secretary, Wilhelm Weierstrass and Theodora born Vonderforst. Karl was the eldest child. He was 12 when his mother died. His father's service was associated with the Tax Department, and therefore, the family had to move from one place to another quite often. His father was a genteel person. He taught his children French and English. Karl started attending school in Münster, and at the age of 14, he entered the catholic Gymnasium at Theodorianum in Padeborn. In addition to good general education, he obtained a good mathematical training at the Gymnasium: stereometry, trigonometry, Diophantine analysis, series expansion. The schooling was thorough. It was for good reason that when the Franco-Prussian War was over, Bismarck said that the war was won by the schoolteacher. There was a scientific library in the Gymnasium. Weierstrass was known to browse mathematical magazines there, especially Crelle's Journal (Journal für die reine und angewandte Mathematik). Each issue of the Journal incorporated four fascicles. In certain years, even two issues were published. Thanks to this periodicity, authors could discuss certain general topics, which eventuated a dialogue and the atmosphere of cooperation. Over the period of Weierstrass' attendance at the Gymnasium (by 1834), 12 issues of the Journal were published. They contained 30 articles of N. Abel and letters he exchanged with A. Legendre; 34 articles of K. Jacobi; 13 articles H. Gudermann, future teacher of Weierstrass. They were essentially devoted to the theory of elliptic functions, which determined the area of Weierstrass' academic interest for life. As he admitted afterwards, he was just carried away by the elliptic functions and creative process in works of Abel, Jacobi and Gudermann.

In addition to these authors, in those years, Crelles Journal published articles of K. Gauss, P. Lejeune-Dirichlet, G. Liouville, A. Legendre, E. Kummer, J. Raabe, which had promoted establishment of the German national mathematical school.

***Bonn University***. The family was financially disadvantaged. Karl even had to work alongside his studies helping a tradeswoman selling butter and ham to keep the books of account. At the age of 19, he graduated from school to be rated *primus omnium*, that is to say, the best of all. His father put high hopes on Karl having chosen the carrier of a public official for him, and Karl set off to the University of Bonn to learn cameral sciences, that is, legal, administrative, and economic sciences that were necessary for civil service, although he was not prone to administrative activities. His attendance at the University of Bonn excited Karl only when it concerned wines, duels, and other mischiefs. Karl was a skilful fencer and was very proud throughout his life that he had never been wounded at a duel. He was granted a special rank in

Saxony corps (fraternity), Fuchsmajor (senior freshman). The administrative career did not seem attractive to him. He completed a course of geometry of J. Plücker and cherished memories of the previous lecturer, Professor K.D. von Münchow (1778-1836), mathematician, astronomer, and physicist. I.V. Goethe, a friend of von Münchow, wrote: "Last year, not only did Prof. von Münchow teach our dear princesses[1] mathematics in Jena, he also prepared them for lectures of Professor Weinhardt, observed, and helped them, from time to time visiting them; moreover, he influenced their morale, state of mind, and behaviour; attracted and held their attention, not to mention his other merits in relation to our dear students" [Gabrichevsky, p. 841]. In addition to his proficiency in mathematics and pedagogical excellence, Von Münchow must have possessed superior human qualities, perfect interpersonal skills, and was very kind.

Karl studied at the university but three terms. However, he stayed in Bonn for two years more. Von Münchow encouraged Weierstrass' intention to study math. As Weierstrass himself wrote on 29 February 1840, "the devout wish to get to know these dearest subjects closer had always attracted me to them, and the more I studied them, the more eager I was in my aspiration to try and devote my effort to their study. What is more, I was lucky to see a well-disposed adviser and supervisor in the deceased Professor von Münchow in Bonn. Finally, the ever growing conviction that the choice of my future profession was wrong, as I felt that I had no inclination or capabilities to become a competent cameralist or lawyer, caused me to decide to throw myself into studying something that is in line with my inclinations and from what I can hopefully expect success" [Kochina, 1985, p. 27-28]. Weierstrass studied all by himself; he studied Laplace Celestial Mechanics and Jacobi's *Fundamenta nova theoriae functionum ellipticarum* which addressed the issue of inverse transformation of Abel integrals and integral systems. The difficulties he faced studying them helped him vanquish Gudermann lecture notes on the theory of modular functions one of the students gave to him [Kochina, 1985, p. 24].

Elliptical integrals appeared in geometric and mechanical problems as far back as at the dawn of differential calculus in works of Newton and Leibniz. Mathematicians tried to reduce them to simpler ones. Euler found that they may be added and multiplied as semicircular arcs and logarithms. They were studied by Lagrange, Legendre, Abel, and Jacobi. Weierstrass made it his objective to continue these studies. Many years later, he wrote: "It was of greatest importance for me to learn in my student days of Abel's letter to Legendre published in Crelles Journal. The first mathematical problem I set to myself was the direct development of the form of presentation of the function denoted by Abel as $\lambda(x)$ from a differential equation determining this function; and the successful solution of this problem was determined by the intention to devote my whole life to math; this happened during my seventh (winter of 1837/38) term" [Biermann, 2008].

This was in the letter to Legendre where Abel wrote about the function $y = \lambda(x)$ with the property that $x = \int_0^y \dfrac{dy}{\sqrt{(1 - y^2)(1 - c^2 y^2)}}$. This function could be presented as a quotient of two everywhere convergent infinite series now referred to as theta series: $y = \dfrac{x + a_3 x^3 + a_5 x^5 + a_7 x^7 + ...}{1 + b_2 x^2 + b_4 x^4 + b_6 x^6 + ...}$. There was no demonstration in Abel's letter, and Weierstrass performed it himself. Weierstrass calculated the coefficients of series and applied the same method to other elliptic functions as well.

---

[1] Weimar princesses Maria and Augusta, granddaughters of Paul I.

Those were the disagreements between the desire to engage in mathematics and requirements of his father which led to a severe inner conflict. Over four years, Weierstrass had not passed a single exam and returned back home. He lost weight and looked ghastly. His brother Peter told Mittag-Leffler: "Karl looked awfully bad when he returned back home! It was so painful to see my elder brother in such a state! Four years and no exam!" [Biermann, 2008].

*Münster*. In October 1838, on the advice of one of the family friends, Weierstrass went to Münster Academy where he hoped to take a course of study in a short period of time to obtain a status of a school teacher. Christoph Gudermann, the second greatest (after Jacobi) lecturer in Germany who gave lectures on elliptic functions, was teaching in the Academy. Weierstrass attended only his courses: Analytical Geometry, Infinitesimal Calculus, Modular Functions and Analytical Spherics. Gudermann gave the two latter courses solely for him. This took only one term, and already in autumn 1839, on special permission from Berlin authorities, Weierstrass started reading for state examinations. In spring 1840, Weierstrass was given three tasks: to write a philosophical work in Latin, a mathematical work which included solving of the problems proposed, and a pedagogical composition. Gudermann set three mathematical problems. The first and the main one, "On development of modular functions", met Weierstrass' preferences and was accompanied by a note to the effect that it was quite a complicated problem for a young analyzer in general and was set with the consent of the examination board only on his express request. The second problem was from elementary geometry; the third one, from theory of mechanics. In this work, relying on certain results of Abel and Jacobi, Weierstrass obtained properties of Abelian functions and various expansions thereof, and Jacobi representations.

Gudermann's comments on the way his student solved these problems were as follows: "1° In this work not only did the author satisfy the expectations of the examination board; moreover, by reference to the system of differential equations which were unknown up to date and which will immediately arouse eminent interest of analyzers, and which he derives directly, serially, and partially one after another, he pioneers the way for the theory of modular functions, and following this way, comes up, as may well be expected, not only with familiar representations of these values, but with totally new results as well. *Thereby, the applicant deservedly joins the scientists crowned with glory*[2].

Considering that he hardly knew any of them when he came to the first lecture on modular functions in Münster, his exceptional success in this comparatively new area of analysis is still more surprising. It can be explained not only by science-oriented diligence of the job applicant, but especially by the extraordinary talent he possesses, which in future will undoubtedly contribute to the science unless feathered out. 2° Wholly satisfactory. 3° This work is also satisfactory.

With such exceptional success of the applicant, no further oral testing needs to be done to assess the amount and thoroughness of his mathematical knowledge if he demonstrates that is able to give a lesson in elements of mathematics using well-structured methods. However, it is absolutely undesirable (for himself or for the science) for him to become a teacher of a grammar school. What he really needs, is to be provided with conditions which would enable him to act as an academic associate professor. Gudermann." [Kochina, 1985, p. 29].

As Weierstrass later wrote (in his letter to Schwartz in 1888), if he were aware of these Gudermann's comments, he would have perceived the value of his work and creativeness and would have been more active in his struggle for a position at the university. His work might well catered for a PhD thesis, however, not at that time, when Münster Academy had no graduate school [Elstrodt]. This Weierstrass' work was published but 54 years later.

---

[2] This phrase was not included in the final version of the comment.

The rest of his examination works were rated wholly satisfactory; at test lessons, he demonstrated sufficient knowledge of Latin, Greek and German to teach in junior school but completely flushed lessons in natural sciences (experimental physics, chemistry, mineralogy, botany, and zoology). Such applicant could not become a schoolteacher. In this regard, correspondence with the Ministry started. As a result, he was allowed to teach mathematics and mathematical physics at high school, and Latin, Greek and German, only in junior school. The insufficient knowledge of other subjects was mentioned in his diploma. During the academic year (1841-42) he was on probation (referendary) at Paulinum Grammar School in Munster.

In Münster, he wrote three more works in the theory of functions of a complex variable. In one of them, analytical functions of a variable were determined with the help of algebraic differential equations. In the same 1842, Cauchy proved the existence theorem. However, at that time, Weierstrass knew nothing about it. His work also included other results which were missing in Cauchy's work. This work contained the notion of a uniform convergence and analytic continuation. It was Weierstrass who introduced the definition of an analytical function as uniformly and unconditionally convergent series[3] (Lagrange wrote nothing at all about convergence, while Cauchy and Abel wrote but about unconditional convergence). Weierstrass did not phrase the notion of a uniform convergence; for him, it just arose out of Abel lemma. Weierstrass spoke about *analytical extension of functions* for the first time mentioning that special points might exist with the property that the radius of convergence approaching these points drops to zero. In his third work, Weierstrass obtained expansion of a function in convergent series with negative and positive powers two years before Pierre Alphonse Laurent (1813-1854). Laurent's work was never published. He sent it to the contest of Paris Academy behind time, and it was known only as stated by Cauchy in 1843 [Cauchy]. In this exposition, Cauchy brought back to memory his theorem from the Lecture Notes in Differential Calculus of 1823: "Let us assume that $x$ denotes a real or an imaginary variable; the real or imaginary function of $x$ may be expanded to a convergent series by growing power of this variable, if the module of the variable keeps a value which does not exceed the smallest of the values for which the functions or its derivative ceases to be finite or continuous." Further, Cauchy said that Laurent extended this theorem of his as follows: "Let us assume that $x$ denotes a real or an imaginary variable; the real or imaginary variable $x$ may be presented as a sum of two convergent series, one with integer growing powers of $x$ and the other one, with integer descending powers of $x$; so far, module $x$ takes values (within the interval) at which the function or its derivative remains finite and continuous." In this article, Cauchy assigned to Laurent's theorem the status of an insignificant conclusion of his theorem, although he used this expansion thereafter. According to [Bottazzini, p. 349], thereafter, Weierstrass used to do without this expansion.

These three works of Weierstrass were published in his collected works but more than 50 years after they were written.

***Deutsch-Krone***. In 1842, Weierstrass was appointed assistant teacher of a progymnasium (junior grammar school) in a small town of Deutsch-Krone (now known as Wałcz, Poland). His workload reached 30 hours per week. He had to teach math, physics, German, Botany, History, Geography, Gymnastics, and Calligraphy. It was at the lessons of Calligraphy that the face of letter $p$ appeared – Weierstrass' function which looked as follows: $\wp$. The lessons of Gymnastics were quite a novelty at that time, and teachers had to first learn themselves. For this purpose, in 1844, Weierstrass went to Berlin where he met a geometrician, J. Steiner, and A.L. Crelle (1780-1855), mathematician, architect, founder and editor of the Journal of Abstract and Applied Mathematics (1826). August Crelle was a self-taught person. It was his merit that having created the journal, he thus united German mathematicians. He could distinguish talented authors and published their works. It was he who discovered Abel's talent and engaged him in his journal, published most of Abel's works in his journal, and cared about his destiny. However, Weierstrass' uncertainty whether his own works were worth while prevented him from showing them to Crelle.

---

[3] The first vision of uniform convergence appeared independently in 1847 in works of J. Stocks and F. Seidel. However, they were speaking of arbitrarily slow convergence; the very notion of a uniform convergence had formed by 1870s in works of Heine (1869) and other mathematicians [Medvedev].

The conditions of his life in Deutsch-Krone were depressing – there was no library in the town; with the small salary (348 thalers per year) he earned he could not even buy postage stamps to send his manuscripts to journals. Weierstrass published his first two works in an annual book of reports of Deutsch-Krone progymnasium. Those were *Notes on Analytic Faculties (Factorials)* and *Reduction of a Definite Threefold Iterated Integral*. The first work was associated with research of Crelle whose work contained certain contradictions and who later suggested that Weierstrass should analyze them in another article which was published in 1856 in Crelle's Journal.

Another Weierstrass' article, *About Socratic Method of Study and Applicability thereof in Schooling,* was published in 1844/45 annual report of Deutsch-Krone progymnasium. This article was his graduation thesis in Münster. The Socratic Method used to be called 'maieutics', i.e. obstetric aid. Asking suggestive questions, the teacher would lead the student to a conclusion of his own. This method was opposed by another Greek method, *akroama*, which meant pleasant reading aloud. This method was more often used in lectures for a large audience. Socrates started his studies with one student and brought him into a certain state of mind. Weierstrass wrote: "Socrates could not establish a general method in common for the entire school. However, it would have been great if his spirit from which his influence proceeded conceived the soul of education and instruction everywhere, his acute pursuance of veritas, beauty and goodness, and love of his clear right" [Kochina, 1985, p. 50]. Weierstrass preferred the method of maieutics in his lectures, i.e. involvement of students in research; he demanded intellectual effort, disapproved of the French lecturing method which implied delivering lectures as a final text. This method delivered benefits when he started lecturing in Berlin. These works went unnoticed, as this collection did not come into the view of professionals.

In 1875, Weierstrass recollected his years of teaching in the grammar school as 14 years of exile to the country of Wieleci and Obodryci (Slavic tribes who lived on the south coast of the Baltic sea, the area of Pomerania and Mecklenburg) [Elstrodt, p. 11]. For a long time, he had no scientific contacts at all.

While in Deutsch-Krone, Weierstrass had an incident with an unlucky engagement, where he played the role of a betrayed fiance, the fact whereof was later described by Schwarz [Dugac, 1973b, p. 167]. Weierstrass was sick for a long time. He recovered very slowly, devoting more and more time to research work.

In 1843, the grammar school in Deutsch-Krone was inspected by a senior inspector who spoke highly of Weierstrass in his report. Therefore, Weierstrass' salary was raised a little (to amount to 400 thalers per year), his promotion was submitted followed by transfer to Catholic Gymnasium of Braunsberg (Braniewo, Poland). However, it was five years later that he was actually appointed to this position.

*Braunsberg*. In autumn 1848, Weierstrass started working at Braunsberg (Braniewo) Catholic Gymnasium in East Prussia not far from Königsberg, now the territory of Poland. Conditions were much better there – they had a library and their headmaster encouraged research. Weierstrass worked a lot on scientific articles – mostly at night. One morning he did not appear at the lesson. When the headmaster came to his place, he found him sitting in the light of the lamp preoccupied with work. In 1850, Weierstrass fell badly ill. He could not do any research for two years. He suffered severe headaches and dizziness for 12 subsequent years.

In Braunsberg, Weierstrass wrote *The Contribution to the Theory of Abel Integrals* which was devoted to the inversion problem for the hyperelliptic case. This work was published in 1848/49 annual report of Braunsberg Gymnasium. This work contained a research related to explicit representation of Abelian integrals through theta series of a number of variables. However, these selected works were not noticed either.

When in 1851 Gudermann died, Weierstrass was considered as a nominee to replace him. However, Plücker whose opinion was decisive said: "I don't even know Weierstrass' first name" [Elstrodt]. In fact, Weierstrass was not aware of this lost opportunity. But during his summer holidays, staying at home in Westfalia, he could read Gudermann's opinion about his graduation thesis including the words: "Thereby, the applicant deservedly joins the scientists crowned with glory of researchers". This inspired him to create the work entitled "About the Theory of Abelian Functions" which was written in 1853. In this work, Weierstrass solved the main problem set by Jacobi regarding inverse transformation of Abelian integrals of the first kind. Getting to know Gudermann's opinion, spirited him up. He sent this work to Crelle's Journal where it was published in Volume 47 (1854).

Thanks to this work Weierstrass got recognition. This article attracted attention of mathematicians, was highly praised by Dirichlet, and affected Weierstrass' destiny. Karl Borchardt (1817-1880), associate professor of Berlin University and student of Jacobi, intentionally came to Braunsberg to meet Weierstrass. This was the start of their long friendship. Thereafter, a delegation from Königsberg led by F.J. Richelot (1808-1875), Jacobi's student, visited Braunsberg to award the PhD Diploma honoris causa to Weierstrass. Handing the Diploma to him, Richelot said: "Each of us has found a teacher in Mr. Weierstrass". On his 80[th] birthday, Weierstrass recollected these words as the most precious ones, having noted: "Everything comes to you in this life but too late" [Kochina, 1985, p. 60]. Thanks to this Diploma, Weierstrass was appointed senior teacher in Braunsberg School.

In his article dedicated to Weierstrass' memory, David Hilbert wrote: "The solution to the Jacobi's inversion problem which was for the first time provided by Weierstrass in these works, which was earlier provided by Rieman with respect to any Abelian integrals, and which was thereafter demonstrated by Weierstrass himself in a different way in his lectures, seems to me to be one of the greatest achievement of analysis" [Hilbert, p. 62].

August Leopold Crelle worked for the Ministry of Education as an advisor on mathematical issues. In 1854, in his letter to the Ministry, he wrote about the then just published work of Weierstrass, saying that it was advisable to provide an appropriate position to him. In his second letter of 1855, shortly before his death, Crelle wrote to the Minister that the outstanding talent of Weierstrass needed support. If Weierstrass is not granted a worthy position, "this man already not quite young and already liable to illnesses caused by the double workload of a teacher and researcher would die ahead of time as Abel and Eisenstein. This would have been another distressful loss for the mathematics. There is a lot of outstanding teachers, while real scientists, who are teachers of the science itself, i.e. teachers of teachers, appear but very seldom" [Biermann, 1966, p. 45]. The article of Weierstrass was immediately translated into French and published in 1854 in issue 19 of Liouville Journal.

On 1 February 1855, Weierstrass wrote a letter to the Minister himself, having enclosed the printing copies of his articles and informed the Minister of their approval: "However, the more valuable this approval is for me and the more it encourages me to proceed with double zeal with the completion of the large works I have started, the more sensitive I am to the fact that my weak health status threatens to make it almost impossible if I stay as I currently am" [Biermann, 1966, p. 45]. After a couple of more letters, on 29 September 1855, he was granted a 12 months leave.

**Berlin**. The death of Gauss in 1855 triggered many transfers in German universities. Dirichlet left Berlin to take a position in Göttingen, Kummer left Breslau to take Dirichlet's position in Berlin. Weierstrass was hoping to get a position in Breslau, but Kummer dissuaded him, because he would have to read only classical courses there. On 19 May 1855, Dirichlet wrote a letter to the Minister of Education to petition for Weierstrass. Austria offered personal professorship at any university to Weierstrass and a salary of 2,000 guilders. Weierstrass hesitated. Kummer wrote to Alexander von Humboldt about it and in three days, Weierstrass was offered a position of a professor at the Königlichen Gewerbeinstitut of Berlin and a salary of 1,500 thalers per year (at that time, one thaler equaled 1.5 Austrian guilders; moreover, guilders were rapidly depreciating in Austria). Weierstrass took this position in July 1856. Soon Weierstrass started lecturing in Friedrich-Wilhelms-Universität Berlin, first as an extraordinary professor (upon the petition of Kummer), and was elected to the Königliche Akademie der Wissenschaften. The election to the Academy granted to the professor the right to choose and give courses of lectures based on his own curriculum. Weierstrass settled in Berlin with his two sisters, Clara and Elisa. Two years later, his widowered father moved to him to stay until he died in 1869.

He was 41. Weierstrass gave 12 hours of lectures per week at the Gewerbinstitut and two lectures at the University; engaged in research and publications. In addition, he had certain duties at the Academy and made prepublication reviews at Crelle's Journal. The overstrain manifested itself on 16 December 1861 at the University – Weierstrass fainted during the lecture. He stopped giving lectures at the Gewerbinstitut for 12 months, although remained on its staff till 1864. On 2 July 1864, he became an ordinary university professor instead of the retired Martin Ohm (1792-1872). Weierstrass lectured for 33 years till 1889, whereupon he started preparing his works for publication.

He himself characterized the epoch from 1864 to 1883 as the time of joint effort of Kummer, Kronecker and himself, as an aspiration to enable the youth at the University over two years "to form a general base with a very large spread of the most important mathematical disciplines" [Biermann K.-R. Die Mathematik, p. 123]. This was a "brilliant constellation of three" [ibidem]; Berlin became the centre which attracted the youth from all over the world to study new sections of mathematics. The professor was playing the role of a researcher, in the first place, and thereafter, that of a teacher.

As H. Hankel mentioned in 1869, after Cauchy's death in 1857, "now the principality of mathematics has indisputably moved to Germany, and, although energetic veterans like Chasles and Liouville still exist in France, they have not got sufficient number of worthy followers who would be able to compete with the Germans" [Hankel, p. 29]. What happened in the epoch of Weierstrass, was the creation of a national school with strong leaders and numerous followers.

For 20 years his cooperation with Kummer, Kronecker and Borchardt constituted an amicable alliance (the league of mutual admiration, as they were called), but in 1880s, relations with sensitive and vain Kronecker started to give way, the fact whereof Weierstrass complained of in 1885 in his letter to Kovalevskaya: "what I am missing more and more, is amicable collaboration with colleagues based on harmony in philosophy and sincere mutual recognition. This has been somehow broken in our University for some years already, and I cannot quite understand the reason why. The only thing I know for sure, it is not myself who had caused this.

My friend Kronecker, with whom we used to reach an accord on the most important issues, and Fuchs resist me: the first one, willfully and intentionally, and the other one, in part submitting to the influence of the first one and in part, being insufficiently aware of the importance of the issue concerned. It is not uncommon that I demonstrate some proposition at a lecture which is recognized to be incorrect at another lecture and doesn't stand up to scrutiny" [Weierstrass, 1973, p. 255].

***Arithmetical approach***. It was Gauss who started developing the theory of functions of the 19th century. He knew the entire range of problems, however, he did not publish anything. In 1798 Gauss wrote a work devoted to elliptic functions and kept it home. When in 1827 he look through Jacobi's and Abel's works, he was very surprised. Gauss wrote to Schumacher in 1827: "Jacobi's results constitute part of my own large work I am going to publish some day. This will be an exhaustive work devoted to this issue, provided that the God be willing to make my life longer and bestow strength and peace of mind to me". The second Gauss' letter was to Bessel: "Mr. Abel anticipated many of my thoughts and made my mission easier approximately by one third, having stated the results very rigorously and elegantly. Abel was following the same way as had been in 1798, therefore, there is nothing remarkable in the fact that we obtained such similar results. To my surprise, this similarity is even in the form and sometimes in notations. Therefore, many of his formulas seem to have been copied from mine. However, to avoid misunderstanding, I should add that I cannot remember a single case when I spoke about this research with a stranger" [Gindikin].

By the mid-nineteenth century, Cauchy had developed basic provisions and structure of analysis: the theory of limits, the notion of continuity and convergence[4], enriched the theory of functions with an integral theorem complex variable and theory of residues. He imposed only the condition of differentiability on the analytical function. An arbitrary function could be represented by an integral. Cauchy's works gave rise to two approaches to the development of the theory of functions: geometric approach of Riemann and arithmetical approach of Weierstrass. Riemann's approach enabled you to visualize properties of elliptic functions and isogonal transformations. Weierstrass' approach was analytical in its nature, logically sound, and enabled to rise to higher degrees of abstraction relative to transcendents which were impossible in terms of geometry. His development of the notion of the number, function, continuity, and least upper bound formed basis for further development of the theory. "For him the function is a formal power series, a 'function element' limited by the convergence circle. There is an analytic continuation procedure outside this circle. Hence, everything is based on the theory of series which is, in its turn, based on arithmetical base. This may apply to functions of several variables.

---

[4] Cauchy to a large extent ingeniously stated and summarized the ideas of B. Bolzano [Sinkevich, 2012b].

The method of Riemann is in the first place the method of discoveries; Weierstrass' method is in the first place the method of proof" [Poincare].

Judging from lecture notes, Weierstrass withdrew most results from Abelian identity in his lectures [Abel, p. 54], the fact whereof is confirmed by Tikhomandritsky: "From here he obtains forms of normal integrals of the second and third kinds; ratios similar to those of Legendre in the theory of elliptic functions between period of integrals of the second and third kinds, prime functions and expression of integrals of all three kinds through them, and algebraic functions which depend on the same irrationality; and from here, as a simple consequence, the theorem of Abel. A special case of the latter leads to the solution of Jacobi's problem, that is to say, he expresses through new variables: sums of $\rho$ integrals of the first kind, the sum of integrals of the second and third kind, and based thereon, considers partial derivatives of sums of integrals of the second kind; the latter ones turn to be partial derivatives of a certain auxiliary functions through which everything can be expressed. If we assume this function as a power of number $e$, we will obtain a single-valued finite and continuous function $\rho$ of new variables which possesses properties similar to those of Jacobi's $\Theta$-function. Finally, Weierstrass withdraws its series expansion. Thus, the theory of Abel's transcendences is reduced to the theory of $\Theta$-functions of many variables in a most natural, not artificial, way as other researchers did." [Tikhomandritsky, p. 45].

His lectures and his concept of analytical function aroused great interest throughout the world and initiated quite a number of research efforts. The number of published works in general theory of functions had rapidly grown under the influence of Weierstrass' lectures (although the number of publications devoted to Abelian functions had grown insignificantly).

*Lectures*. The main results of Weierstrass' research were included in his lecture courses which he had never published. According to H.E. Heine, "Mr. Weierstrass' principles are stated directly in his lectures and indirect oral utterances, in scribal copies of his lectures, and are quite widely spread; however, the author's edition thereof has never been published under his control, which bars the perceptual unity"[Heine, p. 172]. Weierstrass believed that scientific knowledge can be transferred only in case of immediate contact with students, where the lecturer must use materials of his own research and the student must be let into the process of search and taught methods of the research. This "individual" method created a strong school, the doctrine of Weierstrass having spread all over Europe.

*Mittag-Leffler's letter*. On 19 February 1875, G. Mittag-Leffler, one of the favourite and most talented Weierstrass' student, wrote to his motherland to a Swedish Professor Holmgren: "I am delighted with my stay in Berlin in terms of science. I have never found more to study in any other place than here. Weierstrass and Kronecker have a feature so unusual in Germany: they tend to avoid as far as possible printed publications. Weierstrass publishes almost nothing, while Kronecker publishes only results without demonstrations.

They state results of their research in lectures. Nowadays mathematics can hardly demonstrate anything which could be compared to the theory of functions of Weierstrass or Kronecker's algebra.

Weierstrass states the theory of functions in a two- or three-year cycle of lectures and builds a complete theory of elliptic functions and applications thereof in Abelian functions, variational calculus, etc. based on the simplest and clearest notions. His system is basically characterized by the fact that it is completely analytical. He seldom uses geometry and if he does, he does so by way of illustration only. To me, this seems to be an undoubted advantage compared to Riemann's or Clebsch' school.

In fact, it is well known that, based on the theory of Riemann's surfaces, a theory of functions may be built absolutely rigorously and that the geometric system of Riemann is sufficient to find properties of Abelian functions which have not been known until now; however, on the other hand, it is insufficient to find properties of higher order transcendent; otherwise, elements of the theory of functions would have been also introduced in the way that is completely alien to them<…>.

Another feature of Weierstrass is that he tends to avoid any general definitions and all demonstrations pertinent to functions at large. For him a function is a formal power series, and he develops everything from it. This, however, seems to be an extremely complicated method, and I am not sure that, generally speaking, this cannot be achieved the way Cauchy and Liouville do, that is based on general and quite rigorous definitions.

Both Weierstrass and Kronecker are notable for a complete clearness and rigidness of demonstrations. At the same time, they inherited the fear of any type of mathematics from Gausse in establishing basic mathematical notions, and this provides simplicity and naturality to their conclusions which were hardly introduced on such regular basis and with such high degree of rigidity before <…>.

At any rate, from the absolutely formal perspective, Weierstrass' method of reading is beneath criticism, and with such a lecture, even the most inconspicuous French mathematician would have been deemed a completely ineffective lecturer. However, if anyone succeeds after a lot of hard efforts to shape Weierstrass' lecture the way he had meant it to be, then everything becomes clear, simple, and consistent. Probably this surprising drawback of his formal talent explains that very few of his numerous students completely understand him and that the kind of literature he develops is still so insignificant. However, this does not prevent him from enjoying almost idolatrous esteem." [Festschrift, p. 213-214].

Gradually, his lectures formed into a cycle of four terms: *Introduction into the Theory of Analytical Functions*, *Introduction into the Theory of Analytical Functions*, *Abelian Functions*, and *Variations Calculus and Applications of elliptic functions*. He read this cycle until the winter term of 1889/90, his last course being the Variations Calculus. The long-term work on the lectures was reflected in lecture notes later published by his students: Killing's notes of 1868, Hurwitz' notes of 1878, and lecture notes of other students. His lectures on Variations Calculus became widely known thanks to Cobb's lecture notes of 1892/93 and Zermelo's thesis of 1894. Some of the lecture courses of Weierstrass were stated by his students, including *Elements of Arithmetic* by E. Kossak (1872) [Kossak] based on materials of lectures of 1865/66)[5], *Lectures on Weierstrass' Theory of Irrational Numbers* by V. Dantscher [Dantscher], *Experience of Introduction into the Theory of Analytical Functions Based on Weierstrass' Principles* by S. Pincherle (based on lecture notes of 1878) [Pincherle], *The Theory of Analytical Functions* by O. Bierman.

C. Caratheodory who made a great contribution into the theory of variational calculus wrote in the German Literary Newspaper in 1928: "For the lifetime of a generation, mathematicians of all countries engaged in variational calculus regretted that the fundamental discoveries made by Weierstrass in variational calculus could not be found in any of his authentic publications. This may have probably been the only case since the onset of book printing that ideas of a great professional who had revolutionized the entire science were brought to the attention of public only through underground channels."[Kochina, 1985, p. 140-141]. Weierstrass' lectures in variational calculus included the theory of absolute and relative maximums and minimums of functions of one or more variables with criteria of differentiation of relative extremums using Lagrange method based on quadric forms.

---

[5] Translated into Russian by I. Krasovsky and published in Kiev in 1885.

The course of 1861 already contained the notion of continuity in the language ε-δ, which was a decisive step in analysis; the notion of neighborhood, the rigorous definition of infinitely small, the definition of a derivative in the form of $f(x+h) - f(x) = f'(x)h + h(h)$, where $(h) = o(h)$. However, at that time he had no theory of irrational numbers. There was but a sketch when Weierstrass said: "But there also exist values which cannot be expressed by a unity or part thereof; to such values, the form of infinite series applies." [Dugac, 1973a, p. 177]. The first couple of lectures were normally devoted to the notion of a number and the four operations with numbers.

The theory of irrational numbers using the limiting point appeared in Weierstrass' works after 1872 when the notion of the limiting point as the point of accumulation appeared in Hankel's (1870) and Cantor's (1872) works as a point in the neighborhood whereof there are infinitely many points of this set). A Weierstrass' limiting point could be already found in the lecture notes of 1874 (lecture notes of G. Hettner, p. 163-170). After Cantor had introduced the notion of an open and closed set, in the course of Weierstrass of 1874, appeared a δ-neighborhood of a point in $R^n$. This led to creation by Weierstrass of his own conception of continuum [Bottazzini, p. 396]. In the same course, Weierstrass introduced the notion of the least upper bound in the theory of irrational numbers. The statement of the theory was enriched year on year, which is evidenced by lecture notes of the following years [Dugac, 1973b]. In details look [Sinkevich, 2014b].

Having built positive rational numbers as collections composed of elements $e_n$ (and multiples thereof) where $ne_n = 1$ and sum $e_n$ equaled the rational number concerned, Weierstrass introduced a definition as follows: number $B$ is a part of number $A$, if each element of $B$ is an element of $A$. On a set of rational numbers, one can determine substitutions which constitute a replacement of population $n$ of individuals of the form $1/n$ with a unity and any element r through $rn \cdot \frac{1}{n}$. A definition of an equality of two numbers can be provided as follows: two numbers $A$ and $B$ are said to be equal, if any part of $A$ can be transformed by way of substituting it into a part of $B$ and vice versa, any part of $B$ can be turned into a part of $A$. In order to define a collection composed of an infinite number of positive rational numbers, a criterion of finiteness is introduced: collection $A$ is called finite if there is a rational positive number $B$ with the property that any finite part of $A$ is contained in $B$. Thereafter, a set of positive rational numbers is supplemented with the above collections $A$, and their equality is defined in the same way as that of collections composed of a finite number of elements, except that the term 'part' in the definition is replaced by the term' finite part'. Weierstrass proved that if a collection has a finite value, it can be decomposed into two; one of them will contain a finite number of elements, the second one will contain an infinite value, and the sum of elements of the second one is less than any preset number [Dugac, 1973a, p. 179-180].

Since 1874, Weierstrass had been developing the notion of an upper bound of a set [Dugac, 1973b, p. 77] created by Bolzano in 1817. At this, Weierstrass used methods of variational calculus [Sinkevich, 2013].

By the time Cantor's works devoted to the theory of sets appeared, Weierstrass had made the core issue of his lectures the notion of continuum as a perfect connected point set for the theory of monogenic functions and procedures for analytic continuation by way of integration of open disks that form an area where power series representing the function converge uniformly. This was closely connected with the concept of analytic continuation; the area was divided into continuums, and if at least one continuum was found with such property, then the function was monogenic. However, if singularities discontinued the analytic continuation, the area could be divided into plenty of continuums along the bounds that formed them. (Cantor understood the continuum as a connected perfect set). In 1883, Mittag-Leffler wrote to Cantor about both concepts: "I quite agree with your definition of a continuum, however, I would like to refer to the fact that Weierstrass calls the continuum a 'quite connected point set'. It will appear from my work that it is sufficient that such quite connected point set has the necessary place of its own in the theory of analytic functions and cannot be replaced by your continuum". [Turner, p. 114]. Mittag-Leffler had analyzed the difference in concepts of Cantor and Weierstrass in 1883 in his letters to Phragmén. In 1884, Phragmén published an article devoted to this issue (Phragmén. En ny sats inom teorien för punktmängder ("A new theorem within the theory of point sets")).

Weierstrass defined the continuity of a function in the neighborhood of a point with the help of ε-δ apparatus he introduced. He grouped properties of continuous functions: 1) If $f(x_0) \neq 0$, then in the neighborhood of $x_0$ one can find such values of $x$ that $f(x)$ will have the same sign as $f(x_0)$. 2) If $y_1 = f(x_1)$, $y_2 = f(x_2)$, and $y_3$ is any number between $y_1$ and $y_2$, then one can find such value of $x = x_3$ for which $y_3 = f(x_3)$, where $x_3$ is between $x_1$ and $x_2$.

Theorems on functions continuous on an interval belong to Weierstrass: 1) A function continuous on an interval $[a, b]$ is limited on this interval. 2) A function continuous on an interval possesses the largest and the smallest value on this interval[6]. 3) The theorem on approximation of a function: for each real function $f(x)$ continuous on interval $[a, b]$, there is a sequence of algebraic polynomials $P_0(x), P_1(x), ..., P_n(x), ...$ which are uniformly convergent on $[a, b]$ to $f(x)$. A constructive proof of this theorem was provided by S.N. Bernstein in 1912.

28 years before Frechet and Hausdorff, the notion of connection and axiomatics of metric and topological space was formed in lectures of Weierstrass [Sinkevich, 2014b]. However, these notions were auxiliary for him, he needed them to develop the idea of analytic continuation and for the purposes of variational calculus. Therefore, they differed from those created by Cantor. The development of these ideas entailed creation of the theory of metric spaces by M. Frechet and F. Hausdorff, and the theory of functionals in works of V. Volterra and G. Ascoli [Kotsier].

The introductory course included the concept of a number and function with the help of power series, continuity and differentiability, analytic continuation, point of singularity, analytical function of several variables, in particular, Weierstrass' preparation theorem on factorization, and contour integrals. Furthermore, together with Kummer, Weierstrass gave research workshops for bright students. In 1872, he gave workshops in Lobachevskian geometry, where Weierstrass introduced his own coordinates. His audience included not only learners from all over Germany, learners from all over Europe came to his workshops. Therefore, his ideas penetrated other countries as well. In 1873/74, he was elected Head of the University.

***Weierstrass as a lector***. The specifics of Weierstrass' teaching method was that in his course, he first gave grounds. He recommended that new students should attend his lectures from the very beginning. His lecturing was not eloquent – his articulation was not perfect, he could mix sheets of his notes, felt ill at ease, could use his umbrella instead of a sponge; he could improvise, make mistakes, prove his theorems again; sometimes he closed his eyes and engrossed in thought. However, he spoke only about his results and his demonstrations. Being an academician, he had the right to give lectures in accordance with his own curriculum using his own results. He proved all theorems included in his course all by himself. Weierstrass stated the theory of elliptic functions at his lectures in two ways: giving one course, he was based on integrals (this was probably the course attended by Mittag-Leffler); next time (and this course was repeated), he would assume the addition theorem as the point of departure. There are many ways to start stating a theory. The one Weierstrass preferred was very interesting; he would ask a question: when does a function admit the addition theorem?

---

[6] This theorem was for the first time stated by Cauchy, while the complete proof was provided by Heine.

This preference can be easily explained: he thus wanted to let students in the area of Abelian functions when he is finished with their theory. He liked this method of presentation of the subject matter due to its commonality which made the dissemination he meant so easy.

According to Schwarz, he demonstrated mathematics as a field of undiscovered problems. His teaching method started to shape in Münster when he took interest in the method of maieutics (Socrates' method); it got stronger and improved at schools where he had to explain various subjects to negligent grammar school students; and it achieved perfection giving lectures in Berlin. Compare: unlike Weierstrass, Cantor did not directly have students, as he solved all his mathematical problems himself.

However, Felix Klein, as a rebellion against the non-geometric approach of Weierstrass refused to attend his lectures, which he thereafter regretted. Klein said that Weierstrass "enjoyed absolute and incontestable authority; all his theories were accepted by his students as unalterable norms of thinking. His intellectual superiority rather suppressed his students than encouraged their own creative work" [Klein, p. 284].

As a rule, first lectures of Weierstrass were attended by many students, up to 250 students, while by the end of the cycle, only 5 or 7 students remained (as compared to a maximum of 13 students registered to attend Riemann's lectures). However, those were learners pretty advanced in math who were capable of independent research. More than 100 of them became University professors.

***Weierstrass' speech***. The method of Weierstrass is expressed in his speech delivered in 1873 when he was taking over the position of the Rector of the University: "The success of academic teaching is based on the teacher continuously inducing his student to engage in independent research. The teacher achieves this by the fact that the very layout of materials when stating the subject and demonstration of guiding ideas shows the student the way which would lead a mature thinker who possesses all observations to new results in the right sequence or to a better substantiation of the already known results.

At this, the teacher will never miss the chance to point to the bounds which have not yet been crossed by the science by that time and to mention those points based on which a further development of science may be expected in the nearest future. The teacher must not withhold his own research in progress, he will let his student into it without hiding his own errors or disappointments he had faced. Frankly speaking, this way, lectures look not very colourful, elegant, but they are more clear for mentally rigid learners (as, for example, topics stated by most French professors in accordance with lithoprinted notes quite adapted to conform with the established program. These notes are even sometimes entrusted to their assistants who are instructed to read them).

In ancient collections of works of scientific institutions which are hardly read, and in extensive scientific letters ancient scientists exchanged, there is a large amount of scientific materials, in which each person who is able to do so, can find many things inducing to conduct a research of his own and at the same time learn a lot of useful things." [Weierstrass, 1873, p. 1327].

***Weierstrass' lecture of 1886***. Let us mention the first lecture Weierstrass gave in the spring term of 1886 as an example. It was entitled "Select Chapters on the Theory of Functions". Weierstrass gave lectures three times per week, each about 60 minutes long, from early May till late July[7]. Students managed to take word-for-word notes of his lecture, which is why

---

we can hear the quoted speech of Weierstrass. These lectures were published comparatively recently, in 1989 [Weierstrass, 1989]. The lecture provided below is devoted to the notion of a function.

"Tuesday, 25. 5. 1886

These lectures were prepared in such manner as to supplement lectures in the theory of analytical functions given in the winter term of 1884/85. The goal meant to be reached was reached, however, using a more synthetic method; certain results were not duly summarized, and the quality of demonstrations was not totally satisfactory. Therefore, it seems useful to speak about various methods underlying the theory of functions after these lectures, observe them historically and critically in order to demonstrate various viewpoints and try to reconcile them; in a word, to demonstrate the historical development trend of math as a science, especially in the sphere of analysis, and thus explain the core notions of the science. Our goal is to show that principles of mathematical science are based on a really strong foundation. However, even introduction into mathematical sciences necessitates study of various problems, which evidences of the importance and validity of this science in the first place. However, the ultimate goal is something that we always have to bear in mind: one of the fundamentals of a science is the aspiration to obtain an assurance of the correctness of the research."

The understanding of a function as an analytic expression differed from that as a mutual dependence of physical quantities. Newton considered temporal variation of a fluxion and fluent, Descartes considered variation of coordinates of a point on a curve to depend on a certain parameter; Leibniz understood the function as an interval connected with a curve (abscissa, ordinate, radius of curvature). Speaking of dynamics of functions, Euler identified the function as an analytic expression."

Mr. Weierstrass wanted to consider historical changes in the notion of a function, identifying two different definitions thereof. The first one belonged to Leibniz, who, studying algebraic curves, calculated via abscissa such arithmetical expressions as ordinate, length of tangent to the crossing point with the abscissa axis. In keeping with that, Jakob Bernoulli and Leibniz started calling the value which can be calculated based on the other one using a certain arithmetical operation or several operations a function.

This led to further development of algebraic equations, functions with rational coefficients of the variable values, while functions of these values were supposed to be still possible, and roots of any algebraic equations were supposed to be calculable with the help of a limited number of arithmetical operations. This was found by Eiler in his Introductio (Introduction into Analysis of Infinitely Small), while in the times of Lagrange and other significant geometricians, they used the definition of a function as a joint law of functional relationship of arithmetical values.

Together with this definition of a function as an arithmetical expression, there was another one closely followed by Johann Bernoulli, Jacob's brother. He pointed out that variables have geometrical and physical origin and no doubt, each value of the variable or several different variables corresponds to a special continuous variation of one variable or one continuous variation which results in its existence; and the mutual dependence of arithmetical expressions clarifies and confirms such existence. In 1837, Lejeune-Dirichlet provided a definition to a numerical function in relation to functions which have a convergent expansion into Fourier series. However, he applied it mainly to sectionally continuous functions.

Mr. Weierstrass wanted to state a general definition of a function which applies to both arithmetical and geometrical (and natural) variables as provided in the works of Carnot, Cauchy, and Dirichlet:

If a variable is related to another one, where the value of one variable corresponds to that of another one, we will say that the latter variable is a function of the first one. For the purposes of clarification, soon we are going to acquaint ourselves with functions which we could determine when an infinitely small variation of an argument corresponded to the same variation of a function. Thus, the function for a singular point has not been determined. As an example, we can provide an express generalization of the notion of a function in a case. We are going to establish a mutual dependence of effects of two planets. Each of them has a definite mathematical centre of gravity, although not constant – no planet is in a steady state. However, each instant it has a well-defined position in the interplanetary medium. The distance between the two planets is measured as a distance between their centres of gravity. It is a variable value, although each instant, this value is sharply defined. Having duly set the origin of coordinates, we can determine the time, having a certain unit of time, so that each instant we will obtain an indication of the time value; each instant corresponding to a sharply defined unit of length, we will determine the distance; and in the same manner, we can express the time value. Thus we have two interdependent values. Each time value $t$ corresponds to value $r$ of the distance between centres of gravity. Now this law is known for each planet relative another one and takes a complete form, so that the location of the centre of gravity in the environment is identified, the fact whereof is agreed upon when the law is still expressed arithmetically, $r$ is in fact a function of $t$. However, without $r$, it is absolutely impossible to determine the very function of $t$. To tell the truth, the law of mutual influence of two planets is unknown for sure; if Newton's hypotheses underlie the reasoning, then the law has been reliably substantiated. But we know it pretty well that the motion is affected by the resistance and force of friction as well. Therefore, our description is approximate. Now, let us raise a question: what assumptions underlie the mutual dependence of $r$ and $t$ if the law of mutual influence of both planets is unknown? This necessitates the existence of a conditional arithmetical expression found somehow, which allows to calculate $r$ for each arbitrary approximation of value $t$. Another evident question is whether the continuity of a function is ensured if the arithmetical correspondence can be given; whether or not it is always possible to ascertain it. This is one question which can hardly be answered in advance. One can even tend to doubt whether there is common sense in the definition of Johann Bernoulli as opposed to a more general definition of Leibniz. Let us show this in terms of comparison, in respect of each relation of both variables of the arithmetical expression respectively, and similarly for both definitions, while the parity remains in force.

Mr. Weierstrass wanted to reconcile the arithmetical and geometrical (physical) vision of a function: We are not intending to provide an arithmetical definition of a function now, we are going to pave the way the other way round; let us assume that there exists a [numerical] value which depends on one or more variables which continuously varies together with them; and let us demonstrate that this definite value will be presented in a certain arithmetical way. This will be done as follows: begin with an unlimited variable $t$ which takes on values from $-\infty$ to $+\infty$, and the same function $r$, so that if we take any arbitrary small value $\delta$, the entire [entire rational] function of $t$ with rational numerical coefficients is determined in such way that the difference between the true value of $r$ and $r$ from the above expression for any value of $t$ will be less than $\delta$. Further we will suppose that it is the <u>function through an infinite series of presentations for individual rational members that is the entire function</u> <u>with rational numerical coefficients</u>, that is to say, the calculation of the arithmetical expression is so rigid that for any preciseness requirement with respect to the value of $t$ the function may be presented with any approximation. Being governed by Bernoulli's reasoning, we <u>content ourselves with a limited preliminary vision of a real variable</u>. We used to suppose that the presentation [of the function] next to Fourier would solve the problem concerned. Meanwhile, it becomes apparent that there is a continuous function which cannot be obtained if defined this way.

A mathematical expression can be found for a rigorously defined continuous function as well. The advantage of this statement is that it indicates the way to develop properties of any function out of basic notions of continuity, as it is important in any research to derive further notions from the basic ones. Having recognized that a function can be presented for each specific case, then this presentation can be found, that is, you must really know that functions allow to be presented analytically. It appears that the function allows presentation, as far as possible, of itself and of the product thereof by certain powers of arguments inside certain limits of integration. The integration itself is a consequence of the supposed continuity."

At his next lecture, on Wednesday 26 May 1886, Weierstrass was demonstrating the following theorem: $\lim_{k=0} F(x,k) = f(x)$, if

$$F(x,k) = \frac{1}{2k\omega} \int_{-\infty}^{+\infty} f(u)\phi\left(\frac{u-x}{k}\right)du \quad \text{and} \quad \omega = \int_{0}^{\infty} \phi(u)du$$ , and thereafter, he passed over to the notion of continuity.

***Publication of works***. Late in 1885, having celebrated his $70^{\text{th}}$ anniversary, Weierstrass asked for a 12 months leave and spent the entire year 1886 together with his sisters in Switzerland. On his return from Switzerland, he started publishing his works. Seven volumes were published in the period from 1894 to 1927. The first three volumes contained both published and not published works of Weierstrass. The fourth one contained lectures in the theory of Abelian transcendents mostly based on the lecture notes made by Hettner and Knoblauch (1875/76). The fifth volume contained lectures in the theory of elliptic functions; the sixth one, lectures on application of elliptic functions. The seventh volume was published in 1927 and contained lectures in variational calculus. In 1988, *Selected Issues in Complex Analysis* were published which contained Weierstrass' lectures of 1886 [Weierstrass, 1989]. In 1975, the earliest Weierstrass' lecture notes were published. He had read these lectures in 1861 in the Industrial Institute, and 18-year-old G. Schwarz had made these notes. They were found by Pierre Dugac in Mittag-Leffler Institute in Sweden [Dugac, 1973b].

***Students***. In 1871, Germany united into a single state, which triggered national reveille which inspired research in mathematics and thereafter in physics. The orientation of scientific activities of Berlin and Goettingen Universities was shaping the research stream.

Weierstrass' first student was Leo Koenigsberger (who was awarded an academic degree in 1860). He continued the research of his teacher in elliptic functions and differential equations. Understanding the conventionality this classification, we call Weierstrass' followers and disciples, who worked in the mainstream of its main areas of research, in chronological order: L. Fuchs, A.N. Korkin, N.V. Bugaev, K.J. Thomae, H.A. Schwarz, M.A. Tikhomandritsky, E. Kossak, V.P. Ermakov, G.M. Mittag-Leffler, E.I. Zolotarev, F.G. Frobenius, L. Gegenbauer, S.V. Kovalevskaya, F. Schottky, A.V. Vasiliev, K. Runge, O. Bolza, P.M. Pokrovsky, A. Hurwitz, O. Holder, M. Lerch, A. Kneser.

In other areas worked P.G. Bachmann, E. Lampe, F. Mertens, J. Lurot, W. Killing, A.M. Schoenflies, D.F. Selivanov, K. Runge,

The creators of their own trends: M.S. Lie, G. Cantor, F. Klein, E. Husserl, H. Minkowski.

The influence of Weierstrass reached Hermite's students: A. Poincare, G. Darboux, E. Picard, E. Goursat. In Italy, his ideas were followed by S. Pincherle, F. Kasorati, F. Briochi, U. Dini [Sinkevich, 2012a], and J. Peano [Borgato].

***Sofia Kovalevskaya***. Sofia Kovalevskaya (1850-1891) was Weierstrass' favorite student. Having come to him in 1870, she persuaded him to tutorize her as she was not allowed to attend lectures at the University. Having made sure she was a clever and well-trained girl (she had attended Koenigsberger's course of lectures in elliptic functions in Goettingen and had solved some Weierstrass tasks), Weierstrass started with lectures in hyperelliptic functions. She came to him twice a week, and he came to her once a week. In 1872, he taught her variational calculus. Her gratifying attention inspired new mathematical thoughts in him. He called this student his only real friend and considered himself her shepherd. Since 1884, she was teaching at the University of Stockholm. Her success in the research of rigid rotation about a fixed point dates back to 1886. This success was rewarded in 1888 when she was awarded the prize of Paris Academy of Sciences. During her winter holidays of 1890/91, Kovalevskaya was to Berlin. Having returned to Stockholm, she caught cold, went sick, and on 10 February 1891 died at the age of 41. Weierstrass was so shaken by the death of his favorite student that his family was afraid for his own life. He sent a wreath of white lilies to her funeral with an inscription on the mourning ribbon: "To Sonia from Weierstrass". He burnt letters of Kovalevskaya, but his letters to her survived and were published [Kochina, 1981; Weierstrass, 1973].

***Hermite***. Charles Hermite was the leading mathematician of France and considered himself to be Weierstrass' student. He wrote about it to Kovalevskaya on 27 January 1882: "We were taught by the same teacher. It was Mr. Weierstrass. And the main goal of our lectures in Sorbonne and in École Politechnique is to give an account of his works and his great advances to our students. Moreover, you, my gracious lady, form a link in affection between myself and the great geometrician" [Hermite, p. 654].

***Mittag-Leffler***. Magnus Gösta Mittag-Leffler was one the brightest Weierstrass' students (1846-1927). Having graduated from the University in Uppsala in 1873-76, he went abroad to improve his mathematics. In Paris, Hermite advised him to go to Weierstrass, and in 1874/75, he became his student. Mittag-Leffler called Weierstrass 'his great teacher and paternalistic friend'. [Turner, p. 52].

Weierstrass wrote to Kovalevskaya on 15 August 1878: "Mittag-Leffler was a very pleasant student of mine; in addition to thorough knowledge, he possessed a portentous ability to assimilate the subject and a mind focused in the ideal. I am sure that association with him would have a stimulating effect on you". [Weierstrass, 1973, p. 218]. In the same letter, Weierstrass spoke about Mittag-Leffler's situation in Helsingfors University: "They go there further that wherever else in creating *Finnish national* mathematics. And whereas each term in the period of Leffler's stay there local newspapers publish lead articles against Weierstrass', Leffler acts careless mentioning my name in his lectures and articles more often than necessary". [Weierstrass, 1973, p. 218]. Mittag-Leffler's theorem (1876) appeared as an extension of the problem set and solved by Weierstrass. It got its name in an article of Weierstrass and was mouthed by Hermite when he gave a lecture in Sorbonne [Turner, p. 51]. In his letter of 16 December 1874, Weierstrass wrote to Kovalevskaya that with regard to his lectures he was thinking of one unsolved problem: "If an arbitrary infinite series of numbers $a_1, a_2, ..., \infty$, is taken, then the question is: will such entire transcendental function of one variable $x$ always exist with such property that with $x = a_1, a_2, \ldots$ it will disappear, while at any other value it will not? <…> To give an affirmative answer to this question, the condition that as soon as $n$ exceeds a certain limit, the absolute value of $a_n$ will be larger than an arbitrarily assigned value turns to be necessary." [Turner, p. 51].

Weierstrass proved the sufficiency of this condition as well, having presented the required function as

$$\prod_{v_n} E\left(\frac{x}{a_{v_n}}\right)_{v_n}, \text{ where } v_n \text{ is a whole positive number, in particular, } v_n = n, \text{ and}$$

$$E(x)_{v_n} = (1-x)\exp\left(x + \frac{x^2}{2} + \ldots + \frac{x^{v_n}}{v_n}\right).$$

Those were *primary factors*. Poincare believed that their discovery was Weierstrass' main contribution to the theory of functions. Weierstrass' article entitled "Toward the theory of single-valued analytical functions" with this and other results was published in 1876.

Weierstrass had stated the theorems put forward in this article as far back as in summer giving lectures in *Introduction in the Theory of Analytical Functions*. Mittag-Leffler was among the learners, and these lectures inspired him to set a similar problem in the case that for a rational function, instead of zeros, "constants of infinity points" were assumed (the main parts). In 1876, he published two communications containing the so-called Mittag-Leffler's theorem on extension of a meromorphic function: "For any array of numbers $\beta_n$ $(n=1,2,\ldots)$ belonging to a plane of complex numbers which has no points of accumulation in it, there exists a meromorphic function $G$ with poles in points $\beta_n$ and only in these points, and main parts of this function in points $\beta_n$ are the same as the predetermined polynomials of $1/(z-\beta_n)$. In this event, function $G$ can be presented in general terms as an infinite sum of meromorphic functions, each of which has a pole in one point only".

***Summary***. Weierstrass deserves the credit for the creation of a strongly valid mathematical analysis and the theory of elliptic and Abelian functions, and variational calculus. Weierstrass developed the theory of entire and meromorphic functions, provided a canonical presentation of the entire function which has a finite or infinite number of zeros. Instead of three functions of Jacobi, his system of elliptic functions had only one function $\wp(u)$, the simplest one. Weierstrass identified the essential features of algebraic curves which do not change in the case of birational transformations and which are now called 'Weierstrass' points'. Weierstrass developed not only the theory of hyperelliptic integrals, he also studied general Abelian integrals which depend on irrationality: "If an irreducible algebraic equation between two variables allows for an infinite series (eine Schaar) of rationally and uniquely reversible transformations into themselves, then the rank of the algebraic transform is a zero or unity" (a letter to Schwarz).

In 1876, in his article entitled "The theory of single-valued analytical functions" [Weierstrass, 1876], Weierstrass proved the following theorem: if $f(z)$ has the nature of an entire rational function in the neighborhood of each finite point, then it can be presented as a quotient of two entire functions. In the same article, he introduced primary factors and the following theorem: in the neighborhood of an essentially singular point $c$, function $f(x)$ can approximate to any predetermined number arbitrarily closely; at $x = c$, it has no defined value.

(We use to call it the theorem of Sokhotsky-Weierstrass, as eight years earlier, this theorem had been provided by F. Kasorati and Y.V. Sokhotsky independently of each other [Ermolaeva]).

Weierstrass demonstrated that it was possible to construct a single-valued function based on these zeros and a single-valued function with this number of singular points.

The research of Weierstrass extended to the case of multivariable functions. We mean Weierstrass' *preparation theorem* stated in 1886 in *Essays on the Theory of Functions*: Assume $F(x, x_1, x_2, ..., x_n)$ is an analytical function in the neighbourhood of the initial point; let $F(0, 0, ..., 0) = 0$, $F_0(x) = F(x, 0, ..., 0) \neq 0$ and let $p$ be such whole number that $F_0(x) = x^p G(x)$, $G(0) \neq 0$. Then there exists a "selected" polynomial $f(x, x_1, ..., x_n) = x^p + a_1 x^{p-1} + ... + a_p$ whose coefficients are analytical functions $a_j(x_1, ..., x_n)$ in the neighbourhood of the initial point, and function $g(x, x_1, ..., x_n)$, which is an analytical nonzero function in the neighbourhood of the initial point with such property that $F = f \cdot g$ in the neighbourhood of the initial point. It follows from the preparation theorem that if $n > 1$, unlike the case with one complex variable, in any neighborhood of any zero analytical function, there is an infinitude of its zeros. Weierstrass used to include this theorem in his lectures since 1860; it was presented in the lithoprinted publication of 1879.

The theory of Abelian functions was not completed by Weierstrass in whole. Weierstrass introduced the notion of *Abelian* functions, i.e. 2$p$-periodic meromorphic functions $p$ of variables, based on Jacobi's inversion theorem. In 1869, Weierstrass stated a fundamental theorem that there is an algebraic relation between $p+1$ Abelian functions with similar periods; however, he failed to provide a proof [Weierstrass, 1869]. In the following decades, he returned back to this theorem but in vain, because the presentation of meromorphic functions became more complicated with growing dimensionality. Now this problem has been solved [Festschrift, p. 123].

On 18 July 1872, Weierstrass pointed out examples of continuous functions of a real variable which had no definite derivative regardless of the value of this variable (Weierstrass' function: $w(x) = \sum_{n=0}^{\infty} b^n \cos(a^n \pi x)$, where $a$ is an arbitrary odd number which does not equal a unity, and $b$ is a positive number which is less than a unity. This function was created as a contrary instance of Ampere hypothesis).

In 1880, he demonstrated in his work entitled "Zur Functionentheorie" that it was possible to construct such convergent series which would present different functions in different areas. Those were the series which led him to continuous functions having no derivative at all.

Dirichlet's principle got this name in 1851 in the doctoral dissertation of Riemann, Dirichlet's student. Dirichlet implicitly used the principle of existence of a minimum in his lectures and never proved it. Weierstrass showed that in certain situations, this principle was wrong. According to Weierstrass, the assumption that among admissible functions there has to exist one at which the integral must possess the least value has not been proved in terms of mathematics. Based on the assumption that electric current spreads in a conductor, Riemann believed that a problem which is "reasonable in terms of physics", will be "reasonable in terms of mathematics" as well. In 1869, Weierstrass constructed a famous contrary instance. In 1889, Cesare Arzela picked up his idea.

Transcendence of number $e$. In 1882, F. Lindeman proved that number $e^\alpha$ is transcendent for any nonzero algebraic $\alpha$, and in 1885, Weierstrass proved a more general statement currently known as Lindeman-Weierstrass' theorem.

Traditions of Weierstrass' school were fruitful. The doctrine of Weierstrass had gained a legislative nature and spread across Europe thanks to his students and followers.

Weierstrass spent the three last years of his life in a wheel-chair; from time to time his servant took him to a park. Surrounded with veneration of his devotees, he died on 19 February 1897 in Berlin.